\documentclass[11pt]{article}
\usepackage{amssymb,amsmath}

\newtheorem{theorem}{Theorem}[section]
\newtheorem{corollary}[theorem]{Corollary}

\newtheorem{lemma}[theorem]{Lemma}
\newtheorem{proposition}[theorem]{Proposition}
\newtheorem{remark}[theorem]{Remark}

\newenvironment{proof}{\begin{trivlist}\item[]{\it
Proof.}}{\hfill$\square$\end{trivlist}}

\def\mc{{\mathbb{C}}}
\def\field{K}
\def\mn{{\mathbb{N}}}
\def\freeinvar{{\mathcal{F}}}
\def\char{{\rm{char}}}
\def\monoms{{\mathcal{M}}}
\def\qmonoms{{\mathcal{M}(q)}}

\def\partitions{{\mathcal{P}}}
\def\tr{{\mathrm{Tr}}}
\def\rep{{\mathrm{rep}}}
\def\parameterid{\langle P\rangle}
\def\mdeg{{\mathrm{mdeg}}}
\def\aug{{\mathcal{A}}}

\begin{document}
\title{Vector invariants of a class of pseudo-reflection groups and 
multisymmetric syzygies\footnote{J. Lie Theory 19 (2009), 507-525.}} 
\author{M. Domokos\thanks{Partially supported by the Bolyai Fellowship, 
OTKA K61116, and the Leverhulme Trust. }
\\ 
\\ 
{\small R\'enyi Institute of Mathematics, Hungarian Academy of 
Sciences,} 
\\ {\small P.O. Box 127, 1364 Budapest, Hungary,} 
{\small E-mail: domokos.matyas@renyi.mta.hu } }
\date{}
\maketitle 
\begin{abstract} 
First and second fundamental theorems are given for polynomial invariants of  a class of pseudo-reflection groups (including the Weyl groups of type $B_n$), under the assumption that the order of the group is invertible in the base field. As a special case, a finite presentation of the algebra of 
multisymmetric polynomials is obtained. Reducedness of the invariant commuting scheme is proved as a by-product. The algebra of multisymmetric polynomials over an arbitrary base ring is revisited.
\end{abstract}


\section{Introduction}\label{sec:intro} 

Fix natural numbers $n$ and $q$, and a field $\field$. 
Apart from Theorem~\ref{thm:multilinear} and Section~\ref{sec:arbitrarybase}, 
we shall assume that $n!q$ is invertible in $\field$, and assume that $\field$ contains a primitive $q$th root of $1$. 
Denote by $G=G(n,q)$ the subgroup of $GL(n,\field)$ consisting of the monomial matrices whose non-zero entries are $q$th roots of $1$. 
The order of $G$ is $n!q^n$, and as an abstract group, $G$ is isomorphic to the wreath product of the cyclic group $C_q$ of order $q$ and the symmetric group $S_n$; that is, 
$G$ is isomorphic to a semi-direct product 
$(C_q\times\cdots\times C_q)\rtimes S_n$. 

Consider the natural action of $G$ on $V=\field^n$. 
Since $G$ is generated by pseudo-reflections, by the Shephard-Todd Theorem 
\cite{shephard-todd} (see \cite{chevalley} for a uniform proof in characteristic zero, and \cite {smith} for the case when $\char(\field)$ is positive and 
co-prime to the order of $G$) 
the algebra $\field[V]^G$ of polynomial invariants is generated by algebraically independent elements. 
Now consider the diagonal action of $G$ on $V^m=V\oplus \cdots\oplus V$, the direct sum of $m$ copies of $V$. The algebra $\field[V^m]^G$ is no longer a polynomial ring if $m\geq 2$ (and $G$ is not the trivial group). 
In the present paper we show a very short and simple argument that yields  simultaneously the generators of $\field[V^m]^G$ (first fundamental theorem) and the relations among these generators (second fundamental theorem). 
Our main result is Theorem~\ref{thm:finitepresentation}, which 
provides an explicit finite presentation of $\field[V^m]^{G(n,q)}$ in terms 
of generators and relations. 
In the proof we apply Derksen's degree bound on syzygies \cite{derksen} 
and ideas of Wallach and Garsia \cite{haiman}. 

In the special case $q=1$ we have $G=S_n$, and 
$\field[V^m]^{S_n}$ is the algebra of multisymmetric functions, 
which received much attention in the literature (see the references in Remark~\ref{rem:veronese} (ii)). Our approach gives new insight even in this special case, especially by its simplicity and transparency.   We mention also that no finite presentation of the algebra of multisymmetric functions appeared in prior work (apart from the case of $\char(\field)=2$ studied in \cite{feshbach}). 

Another interesting special case is when $q=2$, and the group $G$ is the Weyl group of type $B_n$. The generators of $\field[V^m]^{G(n,2)}$ were determined in 
\cite{jerjen} and \cite{hunziker}. To the best of our knowledge, the second fundamental theorem 
has never been considered in the literature when $q>1$.

The fundamental relation appearing here can be deduced from the theory of trace identities of matrices. This observation leads to the corollary that the 
$GL(n,\mc)$-invariant commuting scheme is reduced, see 
Theorem~\ref{thm:reduced}. 

The present paper joins the content of our preprints \cite{d:msym} and 
\cite{d:bnrel} (some digressions from the preprints have been omitted). 
In addition, in Section~\ref{sec:arbitrarybase} we adjust 
our method for arbitrary base rings and clarify and strengthen the 
known results in this case. In particular, in Theorem~\ref{thm:sigmarelations} we give a new characteristic free (infinite) presentation of the ring of multisymmetric polynomials.


\section{An infinite presentation}\label{sec:infinite} 

Denote by $x(i)_j$ the function on $V^m$ mapping an $m$-tuple of vectors  
to the $j$th coordinate of the $i$th vector component. 
The coordinate ring $\field[V^m]$ is the $mn$-variable polynomial ring 
over $K$ with generators $x(i)_j$, $i=1,\ldots,m$, $j=1,\ldots,n$. 

There is a natural multigrading on $\field[V^m]$ preserved by the action of $G$. 
Namely, $f\in \field[V^m]$ has multidegree $\alpha=(\alpha_1,\ldots,\alpha_m)$, if it has degree $\alpha_i$ in the set of variables $x(i)_1,\ldots,x(i)_n$, for all $i=1,\ldots,m$. Write $\field[V^m]^G_{\alpha}$ for the multihomogeneous component of $\field[V^m]^G$ with multidegree $\alpha$. 
Consider the symbols $x(1),\ldots,x(m)$ as commuting variables, and denote by 
$\qmonoms$ the set of non-empty monomials in the variables $x(i)$ whose degree is divisible by $q$. In particular, $\monoms(1)$ is the set of all non-empty monomials in the variables $x(i)$. 
For $w=x(1)^{\alpha_1}\cdots x(m)^{\alpha_m}\in\qmonoms$ set 
$w_{\langle j\rangle}=x(1)_j^{\alpha_1}\cdots x(m)_j^{\alpha_m}$, 
and define 
$$[w]=\sum_{j=1}^n w_{\langle j\rangle}.$$ 
(These are {\it polarizations} of the usual {\it power sum symmetric functions}.) 
Note that the multidegree of $[w]$ is $\alpha$. 
Since for all $j=1,\ldots,n$ and $w\in\qmonoms$, the term $w_{\langle j\rangle}$ is invariant with respect to the normal subgroup $N$ consisting of the diagonal matrices in $G$, the $S_n$-invariant polynomial $[w]$ is $G$-invariant. 
Sometimes we shall think of the $x(i)$ as the generic diagonal matrices  
${\mathrm{diag}}(x(i)_1,\ldots,x(i)_n)$, $i=1,\ldots,m$, and then $[w]$ is nothing but the trace of the diagonal matrix $w$.  

\begin{proposition}\label{prop:basis} 
The products $[w_1]\cdots [w_r]$ with $r\leq n$, $w_i\in \qmonoms$ 
constitute a $\field$-vector space basis of $\field[V^m]^{G(n,q)}$. 
\end{proposition} 

\begin{proof} An arbitrary monomial $u\in \field[V^m]$ 
in the variables $x(i)_j$ can be written as 
$u=u^1_{\langle 1\rangle}\cdots u^n_{\langle n\rangle}$ 
with a unique $n$-tuple $(u^1,\ldots,u^n)$ of monomials in $\monoms(1)\cup\{1\}$. 
The algebra $\field[V^m]^N$ is spanned by the 
$u$ with $u^1,\ldots,u^n\in\qmonoms\cup\{1\}$. 
The action of $S_n$ permutes these monomials, 
and since $G(n,q)$ is generated by $N$ and $S_n$, 
the $S_n$-orbit sums of the $N$-invariant monomials of multidegree $\alpha$ form a basis in $\field[V^m]^G_{\alpha}$. 
For a multiset $\{w_1,\ldots,w_r\}$ with $r\leq n$, $w_i\in \qmonoms$, 
denote  by $O_{\{w_1,\ldots,w_r\}}$ the $S_n$-orbit sum of the monomial 
${w_1}_{\langle 1\rangle}\cdots {w_r}_{\langle r\rangle}$.  
Call $r$ the {\it height} of this {\it monomial multisymmetric function}. 
Set $T_{\{w_1,\ldots,w_r\}}=[w_1]\cdots [w_r]$. 

So the $O_{\{w_1,\ldots,w_r\}}$ with ${\mathrm{multideg}}(w_1\cdots w_r)=\alpha$ 
form a basis in $\field[V^m]^G_{\alpha}$. 
Assume that the multiset $\{w_1,\ldots,w_r\}$ contains $d$ distinct elements with multiplicities $r_1,\ldots,r_d$ (so $r_1+\cdots+r_d=r$), 
then expanding $T_{\{w_1,\ldots,w_r\}}$ as a linear combination of 
monomial multisymmetric functions, the coefficient of 
$O_{\{w_1,\ldots,w_r\}}$ is $r_1!\cdots r_d!$ (which is non-zero by our assumption on the characteristic of $\field$), and all other monomial multisymmetric functions contributing have height strictly less than $r$. 
This clearly shows the claim. 
\end{proof} 

\begin{remark} {\rm After writing \cite{d:msym} we noticed that the special case 
$q=1$ of Proposition~\ref{prop:basis} appears as Corollary 2.3 of \cite{berele} 
(see also Section 3 of \cite{berele-regev} for the interpretation needed here). }
\end{remark} 

Associate with $w\in\qmonoms$ a variable $t(w)$, and consider the polynomial ring 
$$\freeinvar(q)=\field[t(w)\mid w\in\qmonoms]$$ 
in infinitely many variables. 
Denote by 
$$\varphi(q):\freeinvar(q)\to \field[V^m]^{G(n,q)}$$ 
the $\field$-algebra homomorphism induced by the map 
$t(w)\mapsto [w]$. This is a surjection by Proposition~\ref{prop:basis}. 
Next we introduce a uniform set of elements in its kernel. 
For a multiset $\{w_1,\ldots,w_{n+1}\}$ 
of $n+1$ monomials from $\qmonoms$, we associate an element 
$\Psi(w_1\ldots,w_{n+1})\in\freeinvar(q)$ as follows. 
Write $\partitions_{n+1}$ for the set of partitions 
$\lambda=\lambda_1\cup\cdots\cup \lambda_h$ 
of the set $\{1,\ldots,n+1\}$ into the disjoint union of non-empty subsets 
$\lambda_i$, and denote $h(\lambda)=h$ the number of parts of the partition $\lambda$. 
Set 
$$\Psi(w_1,\ldots,w_{n+1})=\sum_{\lambda\in\partitions_{n+1}} (-1)^{h(\lambda)}
\prod_{i=1}^{h(\lambda)}\left((|\lambda_i|-1)!\cdot 
t(\prod_{s\in \lambda_i}w_s)\right).$$ 
(See Section~\ref{sec:calculations} for examples.)

\begin{proposition}\label{prop:fundrel} 
The kernel of $\varphi(q)$ contains the element  
$\Psi(w_1,\ldots,w_{n+1})$ for arbitrary $w_1,\ldots,w_{n+1}\in\qmonoms$. 
\end{proposition}

\begin{proof} One way to see it is to specialize the fundamental trace identity 
of $n\times n$ matrices to diagonal matrices. 
Indeed, let $Y(1),\ldots,Y(n+1)$ be generic $n\times n$ matrices 
(so their entries generate an $(n+1)n^2$-variable commutative polynomial ring). 
For a permutation $\pi\in S_{n+1}$ with cycle decomposition  
$$\pi=(i_1\cdots i_d)\cdots (j_1\ldots j_e)$$ 
set 
$$\tr^{\pi}=\tr(Y(i_1)\cdots Y(i_d))\cdots\tr(Y(j_1)\cdots Y(j_e)).$$ 
Then we have the equality 
\begin{equation}\label{eq:fundtraceid}
\sum_{\pi\in S_{n+1}}{\mathrm{sign}}(\pi)\tr^{\pi}=0 
\end{equation} 
called the {\it fundamental trace identity} of $n\times n$ matrices. 
(This can be obtained by multilinearizing the Cayley-Hamilton identity 
to get an identity multilinear in the $n$ matrix variables $Y(1),\ldots,Y(n)$, and then multiplying by $Y(n+1)$ and taking the trace; see for example \cite{formanek} for details.) 
The substitution $Y(i)\mapsto w_i$ ($i=1,\ldots,n+1$) 
in (\ref{eq:fundtraceid}) yields the result. 
\end{proof}

\begin{remark}\label{rem:fundiagtraceid} 
{\rm We shall benefit in Section~\ref{sec:semisimple} from the fact that 
in the above proof we descend from a statement about non-commuting matrices. 
A more elementary proof is the following: 
express the $(n+1)$th elementary symmetric function in $n+1$ variables in terms of the first $n+1$ power sums using the Newton formulae; then multilinearize this identity, and specialize the $(n+1)$th coordinates to zero.} 
\end{remark} 

\begin{theorem} \label{thm:firstsecond} \begin{itemize}
\item[(i)] The kernel of the $\field$-algebra homomorphism $\varphi(q)$ is the ideal generated by the 
$\Psi(w_1,\ldots,w_{n+1})$, where $w_1,\ldots,w_{n+1}\in\qmonoms$.  
\item[(ii)] The algebra $\field[V^m]^{G(n,q)}$ is 
minimally generated by the $[w]$, where $w\in\qmonoms$ with 
$\deg(w)\leq nq$. 
\end{itemize}
\end{theorem} 

\begin{proof} 
(i) The coefficient in $\Psi(w_1,\ldots,w_{n+1})$ of the term 
$t(w_1)\cdots t(w_{n+1})$ is $(-1)^{n+1}$, and all other terms are products of 
at most $n$ variables $t(u)$. So the relation $\varphi(q)(\Psi(w_1,\ldots,w_{n+1}))=0$ can be used to rewrite 
$[w_1]\cdots[w_{n+1}]$ as a linear combination of products of at most $n$ 
invariants of the form $[u]$ (where $u\in\qmonoms$). 
So these relations are sufficient to rewrite an arbitrary product of the generators $[w]$ in terms of the basis given by Proposition~\ref{prop:basis}. 
This obviously implies our statement about the kernel of $\varphi(q)$. 

(ii) If $w\in\qmonoms$ and $\deg(w)>qn$, then $w$ can be factored as 
$w=w_1\cdots w_{n+1}$ where $w_i\in\qmonoms$. 
The term $t(w)$ appears in $\Psi(w_1,\ldots,w_{n+1})$ with coefficient 
$-(n!)$, which is invertible in $\field$ by our assumption on the characteristic. Therefore the relation 
$\varphi(q)(\Psi(w_1,\ldots,w_{n+1}))=0$ shows that $[w]$ can be expressed as a polynomial of invariants of strictly smaller degree. 
It follows that $\field[V^m]^{G(n,q)}$ is 
generated by the $[w]$, where $w\in\qmonoms$ with $\deg(w)\leq nq$. 
This is a minimal generating system, because if $\deg(w)\leq nq$ for some $w\in\qmonoms$, then $[w]$ can not be expressed by invariants of lower degree, since there is no relation among the generators whose degree is smaller than 
$(n+1)q$ by (i). 
\end{proof} 

\begin{remark}\label{rem:veronese} {\rm (i) In the special case $n=1$ the above theorem gives the (classically known) defining relations of the $q$-fold Veronese embedding of the projective space $\mathbb{P}^{m-1}$ in $\mathbb{P}^d$, where $d=\left(\begin{array}{c}m+q-1\\ q\end{array}\right)-1$. 

(ii) In the special case $q=1$ our group is the symmetric group $S_n$ and 
$\field[V^m]^{S_n}$ is called the {\it algebra of multisymmetric polynomials}. 
It is an old result proved by 
Schl\"afli \cite{schlafli}, 
MacMahon \cite{macmahon}, Weyl \cite{w},   
that when $\field$ has characteristic zero, this algebra is generated by the polarizations of the elementary symmetric polynomials 
(we mention that Noether \cite{noether} used this to prove her general degree bound on generating invariants of finite groups).  
The result on the generators was extended to the case $\char(\field)>n$ by 
Richman \cite{richman} (see also \cite{stepanov}, 
\cite{briand}, \cite{fleischmann2} for other proofs). 
Counterexamples and further results on the generators in the case $0<\char(\field)\leq n$ are due to Fleischmann \cite{fleischmann}, Briand \cite{briand}, Vaccarino \cite{vaccarino}.  
The relations among the generators had been studied classically by Junker  \cite{junker1} \cite{junker2}, \cite{junker3}. 
Working over $\mc$ as a base field, an infinite presentation was given 
by Dalbec \cite{dalbec} (partly reformulating results of Junker \cite{junker2}),  and the special case $q=1$ (and $\field=\mc$) of Theorem~\ref{thm:firstsecond} appears in this explicit form in a recent paper by Bukhshtaber and Rees \cite{bukhshtaber-rees} (using a different language and motivation). 
Working over an arbitrary base ring $\field$,  
\cite{vaccarino} gives $\field$-module generators of the ideal of relations among a reasonable infinite generating system 
(see Section~\ref{sec:arbitrarybase}). 

(iii) In the special case $q=2$ our group $G(n,2)$ is the Weyl group of 
type $B_n$, and the ring of invariants 
$\mc[V^m]^{G(n,2)}$ is isomorphic to the ring of invariants of the special orthogonal group $SO(2n+1,\mc)$ acting via the adjoint representation on the commuting variety of its Lie algebra, 
and is also isomorphic to the ring of invariants of the symplectic group 
$Sp(2n,\mc)$ acting on the commuting variety of its Lie algebra 
(see for example \cite{hunziker} for this connection).   

(iv) In our point of view, Theorem~\ref{thm:firstsecond} (i) is a close relative of the result of Razmyslov \cite{razmyslov} and Procesi \cite{procesi} saying that all trace identities of $n\times n$ matrices are consequences of the fundamental trace identity. Their original approach is based on 
Schur-Weyl duality; Kemer in \cite{kemer} gave an elementary combinatorial 
proof valid on the multilinear level in all characteristic.  

(v) Identify $V$ with the $\field$-algebra of $n\times n$ diagonal matrices. Denote by 
 ${\mathrm{Mor}}_{S_n}(V^m,V)$ the space of $S_n$-equivariant polynomial maps from $V^m$ to $V$. 
 Since $S_n$ acts on $V$ by algebra automorphisms,  ${\mathrm{Mor}}_{S_n}(V^m,V)$ is a $\field$-algebra with pointwise multiplication of maps $V^m\to V$. 
Identify the space ${\mathrm{Mor}}(V^m,V)$ of all polynomial maps from $V^m$ to $V$ with the algebra of $n\times n$ diagonal matrices with entries in $\field[V^m]$ in the obvious way; then we have the usual trace function on  ${\mathrm{Mor}}(V^m,V)$, and it is an algebra with trace in the sense of 
\cite{procesi:1987}. 
The special case $q=1$ of Theorem~\ref{thm:firstsecond} has 
the following corollary (by some methods of \cite{procesi}): if $n!$ is invertible in $\field$, then 
the algebra ${\mathrm{Mor}}_{S_n}(V^m,V)$ of $S_n$-equivariant polynomial maps from $V^m$ to $V$ is the trace-stable subalgebra generated by the generic diagonal matrices $x(1),\ldots,x(m)$ in the algebra of diagonal matrices with entries from $\field[V^m]$. Furthermore, it is 
the free object of rank $m$ 
in the category of commutative algebras with a trace satisfying the $n$th 
Cayley-Hamilton identity (in the sense of \cite{procesi:1987}). 
This latter result is due to Berele, see Theorem 2.1 in \cite{berele}. } 
\end{remark}   

An element of $\field[V^m]$ is {\it multilinear} if it is multihomogeneous with multidegree $(1,\ldots,1)$. Define a multigrading on $\freeinvar(q)$ by setting the multidegree of $t(w)$ equal to the multidegree of $w$. 
The multilinear version of the above theorem is valid over an arbitrary 
base field. Abusing notation we keep on writing 
$\field[V^m]^{G(n,q)}$ even if the group $G(n,q)$ is not defined over $\field$ 
(i.e. when $\field$ does not contain a primitive $q$th root of $1$); in this case the notation refers to the algebra of invariants of $S_n$ acting on the subalgebra of $\field[V^m]$ spanned by monomials 
$u=u^1_{\langle 1\rangle}\cdots u^n_{\langle n\rangle}$, where $u^i\in\qmonoms\cup\{1\}$ for all $i$ (with the notation of the proof of Proposition~\ref{prop:basis}). 

\begin{theorem} \label{thm:multilinear} 
Let $\field$ be an arbitrary field (or even an arbitrary commutative ring). 
Then the multilinear component of $\field[V^m]^{G(n,q)}$ coincides with the multilinear component of the image of $\varphi(q)$. 
Moreover, any multilinear element in the kernel of $\varphi(q)$ 
is contained in the ideal generated by 
$\Psi(w_1,\ldots,w_{n+1})$, where $w_i\in \qmonoms$ and the multidegree of 
$w_1\cdots w_{n+1}$ belongs to $\{0,1\}^m$. 
\end{theorem}  

\begin{proof} The only point in the proof of Proposition~\ref{prop:basis} where 
we need an assumption on $\field$ is to guarantee that $r_1!\cdots r_d!\neq 0$. 
If we consider multilinear elements of $\field[V^m]$ only, then 
all the $r_i$ here are automatically equal to $1$. So the same argument 
yields that 
the multilinear elements in $\field[V^m]$ of the form 
$[w_1]\cdots [w_r]$ (with $r\leq n$, $w_i\in\qmonoms$) constitute a 
$\field$-basis in the multilinear component of $\field[V^m]^{G(n,q)}$. The relations from our statement are obviously sufficient to rewrite an arbitrary multilinear product of elements $[w]$ (with $w\in\qmonoms$) as a linear combination of such products with at most $n$ factors. This implies the claim. 
\end{proof} 

\begin{remark}\label{rem:indecomposable} 
{\rm As an immediate corollary one recovers the 
known fact that the invariant $[x(1)\cdots x(m)]\in\field[V^m]^{S_n}$ is indecomposable (can not be expressed by lower degree multisymmetric polynomials)  when $0<\char(\field)\leq n$. 
Indeed, in this case each term of the fundamental relation is the product of at least two traces, hence there is no relation involving 
$[x(1)\cdots x(m)]$ as a term. 
This has interesting consequences for matrix invariants in positive characteristic, and for vector invariants of the orthogonal group; see \cite{domokos}.  }
\end{remark} 


\section{Finite presentation}\label{sec:finite} 

Building on Derksen's general degree bound for syzygies in \cite{derksen},  
we derive a finite presentation from Theorem~\ref{thm:firstsecond}.  

Denote by $J_m(q)$ the ideal in $\field[V^m]$ generated by the homogeneous 
$G(n,q)$-invariants of positive degree (the so-called Hilbert ideal), and denote by 
$\tau_{G(n,q)}(V^m)$ the minimal natural number $d$ such that all homogeneous elements of degree $\geq d$ belong to $J_m(q)$. 

\begin{lemma}\label{lemma:tau}
We have $\tau_{G(n,q)}(V^m)=1+\sum_{j=1}^n(qj-1)$. 
\end{lemma}

\begin{proof} In the special case $q=1$ (and $\field=\mc$, $m=2$) 
this lemma appears in \cite{haiman} and is due to Wallach and Garsia. 
Define the {\it exponent} of a monomial $\prod_{i,j}x(i)_j^{\alpha_{ij}}$ 
as $(\beta_1,\ldots,\beta_n)$, where $\beta_j=\sum_{i=1}^m\alpha_{ij}$ 
(in other words, $\beta_j$ is the exponent of $x_j$ if we specialize all 
the $x(i)_j$ to $x_j$). Introduce a partial ordering on the set of monomials in 
$\field[V^m]$: a monomial $u$ is smaller than $v$ 
if they have the same multidegree, and the exponent of $u$ is lexicographically smaller than the exponent of $v$. This partial ordering is compatible 
with multiplication of monomials. 
We shall show that if $\beta_k\geq qk$ for some $k\in\{1,\ldots,n\}$, where 
$(\beta_1,\ldots,\beta_n)$ is the exponent of some monomial $u$ in $\field[V^m]$, then modulo $J_m(q)$, the monomial $u$ can be rewritten as a linear combination of smaller monomials. 
We need some details from the proof of the special case $q=1$ of our lemma. 
To simplify notation when $m=1$, write $x_j=x(1)_j$. 
Then it is shown in Lemma 3.2.1 in \cite{haiman} that 
\begin{equation}\label{eq:hk}
h_k(x_k,\ldots,x_n)\in J_1(1) \mbox{ for all }k=1,\ldots,n, 
\end{equation} 
where $h_k$ denotes the $k$th complete symmetric polynomial of the arguments. 
Note that (\ref{eq:hk}) shows that modulo $J_1(1)$, the monomial $x_k^k$ can be rewritten as a linear combination of smaller monomials. 
Applying to (\ref{eq:hk}) the comorphism of the $S_n$-equivariant morphism 
$V^m\to V$, $(v_1,\ldots,v_m)\mapsto v_1+\cdots+v_m$ we get 
\begin{equation} \label{eq:hkpolarize} 
h_k(\sum_{i=1}^mx(i)_k,\ldots,\sum_{i=1}^mx(i)_n)\in J_m(1) 
\mbox{ for all }k=1,\ldots,n. 
\end{equation} 
Since the ideal $J_m(1)$ is multihomogeneous, all multihomogeneous components of the left hand side of (\ref{eq:hkpolarize}) belong to $J_m(1)$. 
In particular, for all $k=1,\ldots,n$, there exists a linear combination $f_k$ of $k$-linear monomials smaller than $x(1)_kx(2)_k\cdots x(k)_k$ such that 
\begin{equation}\label{eq:3} 
x(1)_kx(2)_k\cdots x(k)_k-f_k\in J_k(1).
\end{equation} 
(At this point it is essential that $k!$ is assumed to be invertible in $\field$.) 
Now take an arbitrary $k$-tuple $w^{(1)},\ldots,w^{(k)}\in\qmonoms$. 
The map $x(i)_j\mapsto w^{(i)}_{\langle j\rangle}$ induces a 
$\field$-algebra homomorphism $\gamma:\field[V^k]\to \field[V^m]^N$. 
Clearly $\gamma$ is an $S_n$-equivariant map, whence $\gamma(J_k(1))\subseteq  J_m(q)$. 
Applying $\gamma$ to (\ref{eq:3}) we get that 
\begin{equation}\label{eq:4} 
w^{(1)}_{\langle k\rangle}\cdots w^{(k)}_{\langle k\rangle}-\gamma(f_k)\in J_m(q), 
\end{equation} 
where $\gamma(f_k)$ is a linear combination of monomials smaller than 
$w^{(1)}_{\langle k\rangle}\cdots w^{(k)}_{\langle k\rangle}$.  
Now let $w$ be a monomial in the variables $x(1)_k,x(2)_k,\ldots,x(m)_k$ such that $\deg(w)=qk$. Then $w$ can be factored  as 
$w=w^{(1)}_{\langle k\rangle}\cdots w^{(k)}_{\langle k\rangle}$ with 
$w^{(i)}\in\qmonoms$, whence modulo $J_m(q)$ the monomial $w$ can be rewritten as a linear combination of smaller monomials by (\ref{eq:4}).  
This clearly implies that the factor space $\field[V^m]/J_m(q)$ is spanned by the images of the monomials whose exponent $(\beta_1\ldots,\beta_n)$ satisfies 
$\beta_k\leq qk-1$ for all $k=1,\ldots,n$. Since $J_m(q)$ is a homogeneous ideal, 
the factor space $\field[V^m]/J_m(q)$ inherits the grading from $\field[V^m]$. 
This shows that the homogeneous components of degree $>\sum_{k=1}^n(qk-1)$ are all contained in $J_m(q)$, implying that 
$\tau_G(V^m)\leq 1+\sum_{k=1}^n(qk-1)$.  

In the special case $m=1$, the coordinate ring $\field[V]$ is known to be a free module over $\field[V]^{G(n,q)}$, and the latter is a polynomial algebra in $n$ generators having degrees $q,2q,\ldots,nq$. This implies a formula for the dimensions of 
the  homogeneous components of the graded vector space $\field[V]/J$, called the coinvariant algebra (see for example \cite{chevalley}). In particular, the highest degree non-zero homogeneous component has degree $\sum_{k=1}^n(qk-1)$. 
This shows the reverse inequality  $\tau_{G(n,q)}(V^m)\geq 1+\sum_{k=1}^n(qk-1)$.  
\end{proof} 

For a natural number $d$, consider the finitely generated subalgebra of $\freeinvar(q)$ given by 
$$\freeinvar(q,d)=\field[t(w)\mid w\in\qmonoms,\ \deg(w)\leq d].$$ 

\begin{theorem}\label{thm:finitepresentation}
The kernel of the $\field$-algebra surjection 
$$\freeinvar(q,qn(n+1)-2n+2)\to \field[V^m]^G\mbox{ induced by \quad} 
t(w)\mapsto [w]$$ 
is generated as an ideal by 
the elements $\Psi(w_1,\ldots,w_{n+1})$, where 
$w_1,\ldots,w_{n+1}\in\qmonoms$, and the degree of the product 
$w_1\cdots w_{n+1}$ is not greater than 
$qn(n+1)-2n+2$. 
\end{theorem} 

\begin{proof} A general result of Derksen \cite{derksen} says that 
the ideal of relations in a minimal presentation of $\field[V^m]^G$ is generated in degree $\leq 2\tau_G(V^m)$. 
Therefore our statement follows from Theorem~\ref{thm:firstsecond},  
Lemma~\ref{lemma:tau}, and the general Lemma~\ref{lemma:redundant} below. 
\end{proof} 

For a finitely generated commutative graded $\field$-algebra $A$ denote 
by $\tau(A)$ the minimal non-negative integer $\tau$ such that $A$ and its 
first syzygy ideal are generated in degree $\tau$.     
That is, let $\rho:\field[x_1,\ldots,x_r]\to A$ be a surjective $\field$-algebra homomorphism such that 
$\rho(x_i)$, $i=1,\ldots,r$, is a minimal homogeneous generating system of $A$, and the grading on the polynomial algebra $\field[x_1,\ldots,x_r]$ is defined so that $\rho$ preserves the grading. Let $f_1,\ldots,f_s$ be a minimal homogeneous generating system of the ideal $\ker(\rho)$. If $\ker(\rho)=0$, then define 
$\tau(A)=\max\{\deg(x_i) \mid i=1,\ldots,r\}$, and otherwise set 
$\tau(A)=\max\{\deg(f_j)\mid j=1,\ldots,s\}$.  

\begin{lemma}\label{lemma:redundant} 
Let $\varphi:\field[t_{\mu}\mid \mu \in M]\to A$ be a surjective homomorphism of graded algebras from a not necessarily finitely generated polynomial algebra, and let 
$\{f_{\beta} \mid \beta\in B\}$ be a homogeneous generating system of the ideal 
$\ker(\varphi)$. Then the kernel of the restriction of $\varphi$ to the subalgebra 
$\field[t_{\mu}\mid \deg(t_{\mu})\leq \tau(A)]$ is generated as an ideal 
by $\{f_{\beta}\mid \deg(f_{\beta})\leq \tau(A)\}$. 
\end{lemma} 

\begin{proof} 
Denote by $I$ the ideal of $R_1=\field[t_{\mu}\mid \deg(t_{\mu})\leq \tau(A)]$ 
generated by $\{f_{\beta}\mid \deg(f_{\beta})\leq \tau(A)\}$. 
Take a subset $C\subset M$ such that the restriction 
$\varphi_0:R_0=\field[t_{\mu} \mid \mu\in C]\to A$ of $\varphi$ 
is a minimal presentation of $A$. Then $\ker(\varphi_0)$ is generated in degree 
$\tau(A)$ and is contained in $\ker(\varphi)=(f_{\beta} \mid \beta\in B)$, 
hence $\ker(\varphi_0)\subset I$. Moreover, since 
${\mathrm{im}}(\varphi_0)={\mathrm{im}}(\varphi)$, 
for each $\mu\in M$ there is an element $g_{\mu}\in R_0$ with $\deg(t_{\mu})=\deg(g_{\mu})$ such that 
$t_{\mu}-g_{\mu}\in \ker(\varphi)$. Obviously, if $\deg(t_{\mu})\leq \tau(A)$, 
then $t_{\mu}-g_{\mu}\in I$. It follows that $R_1=R_0+I$, implying 
$\ker(\varphi_1)=\ker(\varphi_0)+I=I$. 
\end{proof} 

\begin{remark} {\rm Since the Hilbert ideal of a subgroup $H$ of $G(n,q)$ contains the Hilbert ideal of $G(n,q)$, we have $\tau_H(V^m)\leq \tau_{G(n,q)}(V^m)$. 
Thus Lemma~\ref{lemma:tau} provides upper bounds on the degrees of the relations 
in a minimal presentation of $\field[V^m]^H$ by the result of Derksen \cite{derksen} cited above. Among subgroups of $G(n,q)$ are other series of pseudo-reflection groups, including the Weyl groups of type $D_n$ or the dihedral groups. }
\end{remark}


\section{Semisimple commutative representations}\label{sec:semisimple} 

In this section we work over the field $\mc$ of complex numbers 
(what is essential for Theorem~\ref{thm:reduced} is that the base field has characteristic zero). 
Write $F_m$ for the free associative algebra $\mc\langle x_1,\ldots,x_m\rangle$, 
and $A_m$ for the commutative polynomial algebra $\mc[x(1),\ldots,x(m)]$. 
Denote by $\rep(F_m,n)$ the space of $m$-tuples of complex $n\times n$ matrices, endowed 
with the simultaneous conjugation action of $GL(n,\mc)$.  
The points of this affine variety determine $n$-dimensional representations of $F_m$ in the obvious way, and the $GL(n,\mc)$-orbits are in a one-to-one correspondence with the isomorphism classes of $n$-dimensional representations of $F_m$. The algebraic quotient 
$\rep(F_m,n)/\!/GL(n,\mc)$ parameterizes the isomorphism classes of semisimple 
$n$-dimensional representations of the free algebra, see \cite{artin}. 
The coordinate ring $\mc[\rep(F_m,n)]$ is the $n^2m$-variable polynomial algebra generated by the entries of the generic $n\times n$ matrices 
$Y(1),\ldots,Y(m)$. 
Denote by $J$ the ideal in $\mc[\rep(F_m,n)]$ generated by the entries 
of the commutators $Y(i)Y(j)-Y(j)Y(i)$, $1\leq i<j\leq m$. 
The quotient algebra $\mc[\rep(F_m,n)]/J=\mc[\rep(A_m,n)]$ is the coordinate ring of the {\it scheme $\rep(A_m,n)$ of $n$-dimensional representations of 
$A_m$} (by definition of this affine scheme). 
It is a long standing open problem in commutative algebra 
whether this scheme is reduced or not, or equivalently, whether $J$ is a radical ideal or not; see for example \cite{hreinsdottir}. 
The common zero locus of $J$ (i.e. the set of $\mc$-points of the scheme 
$\rep(A_m,n)$) is the so-called {\it commuting variety} consisting of $m$-tuples of pairwise commuting $n\times n$-matrices. So the question is whether $J$ is the whole vanishing ideal of the commuting variety or not.  
The embedding of the space $V^m$ of $m$-tuples of diagonal matrices into the commuting variety induces a surjective homomorphism 
$\beta:\mc[\rep(A_m,n)]\to\mc[V^m]$. 
Consider the homomorphism 
$\gamma:\freeinvar(1)\to \mc[\rep(A_m,n)]^{GL(n,\mc)}$ given by 
$t(w)\mapsto \tr(Y(1)^{\alpha_1}\cdots Y(m)^{\alpha_m})$ 
for $w=x(1)^{\alpha_1}\cdots x(m)^{\alpha_m}\in\monoms(1)$, 
where we keep the notation $Y(j)$ for the images of the generic matrices under the natural surjection 
$$M(n,\mc[\rep(F_m,n)])\to M(n,\mc[\rep(A_m,n)].$$  
It is well known that $\gamma$ is surjective, see \cite{sibirskii}; 
we use here that any $GL(n,\mc)$-invariant in $\mc[\rep(A_m,n)]$ lifts to 
an invariant on $\rep(F_m,n)$,  
since our base field has characteristic zero. 
The fundamental trace identity holds in $M(n,C)$ for an arbitrary commutative ring $C$ (this is equivalent to the Cayley-Hamilton theorem), 
therefore $\Psi(w_1,\ldots,w_{n+1})\in\ker(\gamma)$ for any 
$w_i\in\monoms(1)$. 
It follows by Theorem~\ref{thm:firstsecond} that $\gamma$ factors through a 
surjection 
$$\overline{\gamma}:\mc[V^m]^{S_n}\to\mc[\rep(A_m,n)]^{GL(n,\mc)}, \quad   
[w]\mapsto \tr(Y(1)^{\alpha_1}\cdots Y(m)^{\alpha_m})$$  
for $w=x(1)^{\alpha_1}\cdots x(m)^{\alpha_m}\in\monoms(1)$. 
So we have the surjections 
$$\mc[V^m]^{S_n}\stackrel{\overline{\gamma}}\longrightarrow
\mc[\rep(A_m,n)]^{GL(n,\mc)}\stackrel{\beta}\longrightarrow
\mc[V^m]^{S_n}.$$
Since $\beta\circ\overline{\gamma}$ is the identity map of 
$\mc[V^m]^{S_n}$ by definition of the  maps $\overline{\gamma}$ and $\beta$, 
we conclude that they are isomorphisms. Thus we obtained the following result, supporting the conjecture that the commuting scheme is reduced: 

\begin{theorem}\label{thm:reduced} 
The surjection $\mc[\rep(A_m,n)]\to\mc[V^m]$ induced by the embedding of 
the space of $m$-tuples of diagonal matrices into the commuting variety 
restricts to an isomorphism 
$\mc[\rep(A_m,n)]^{GL(n,\mc)}\cong 
\mc[V^m]^{S_n}$. 
In particular, the radical of  $\mc[\rep(A_m,n)]$  
contains no non-zero $GL(n,\mc)$-invariants. 
\end{theorem} 

\begin{remark}  {\rm 
One may paraphrase the above theorem by saying that 
"the scheme of $n$-dimensional 
semisimple representations of the commutative 
$m$-variable polynomial algebra is reduced". 
(The weaker statement that the algebraic quotient of the commuting variety 
with respect to the action of $GL(n,\mc)$ is isomorphic to $V^m/S_n$ is well known, and is explained for example in Proposition 6.2.1 of \cite{haiman}.) 
In the special case $m=2$, reducedness of  $\mc[\rep(A_m,n)]^{GL(n,\mc)}$ was proved by different methods by Gan and Ginzburg, see Theorem 1.2.1 in \cite{gan-ginzburg}. 
For arbitrary $m$ it appeared in our preprint \cite{d:msym}, 
and is proved also by Vaccarino in his parallel preprint [arXiv:math.AG/0602660], that appeared in the meantime as \cite{vaccarino2}. }
\end{remark}


\section{Multisymmetric functions over an arbitrary base ring} 
\label{sec:arbitrarybase}

Throughout this section $\field$ is an arbitrary base ring. 
We extend the method of Section~\ref{sec:infinite} to this case. 

For a monomial $w\in\monoms=\monoms(1)$ and $l\in\{1,\ldots,n\}$ 
we set 
$$\sigma_l(w)=\sum_{1\leq i_1<\cdots<i_l\leq n}
w_{\langle i_1\rangle} \cdots w_{\langle i_l\rangle}.$$ 
In other words, $\sigma_l(w)$ is the $l$th coefficient of the characteristic polynomial of the diagonal matrix $w$. 
By the multiset $\{w_1,\ldots,w_d;r_1,\ldots,r_d\}$ of monomials from $\monoms$ 
we mean the multiset which contains $d$ distinct monomials $w_1,\ldots,w_d$, and the multiplicity of $w_i$ is $r_i$. 

\begin{proposition}\label{prop:basis2} 
A free $\field$-module basis in $\field[V^m]^{S_n}$ is formed by 
the products 
$\sigma_{r_1}(w_1)\cdots\sigma_{r_d}(w_d)$,  
where $\{w_1,\ldots,w_d;r_1,\ldots,r_d\}$ 
range over all multisets of monomials from $\monoms$ 
with $r_1+\cdots+r_d\leq n$.  
\end{proposition} 

\begin{proof} The proof is essentially the same as that of Proposition~\ref{prop:basis}: when we expand 
$\sigma_{r_1}(w_1)\cdots\sigma_{r_d}(w_d)$ as a linear combination of monomial multisymmetric functions, the coefficient of 
$O_{\{w_1,\ldots,w_d;r_1,\ldots,r_d\}}$ is $1$, and all other monomial multisymmetric functions contributing have strictly smaller height. 
This implies the claim, since the $O_{\{w_1,\ldots,w_d;r_1,\ldots,r_d\}}$ with 
$r_1+\cdots+r_d\leq n$ form a $\field$-basis in $\field[V^m]^{S_n}$. 
\end{proof} 

In particular, the $\field$-algebra $\field[V^m]^{S_n}$ is generated by 
the $\sigma_r(w)$. 
Our aim is to determine generators of the ideal of relations among these generators. With $w\in\monoms$ and $r\in\{1,\ldots,n\}$ associate a commuting 
indeterminate $e_r(w)$, and denote by 
$$\phi:\field[e_r(w)\mid w\in\monoms,r=1,\ldots,n]\to \field[V^m]^{S_n}$$
the $\field$-algebra surjection induced by 
$e_r(w)\mapsto \sigma_r(w)$. 

Denote by $\sigma_l^{\alpha}$ the multihomogeneous component of multidegree $\alpha=(\alpha_1,\ldots,\alpha_s)$ of the $l$th characteristic coefficient of the diagonal matrix $x(1)+\cdots +x(s)$ (note that necessarily $l=\alpha_1+\cdots+\alpha_s$). Specializing Amitsur's formula \cite{amitsur} for the characteristic coefficients of a linear combination of matrices one obtains 
formula (\ref{eq:amitsur}) below. To state it we need to introduce some 
integers. Take non-commuting variables $Y_1,\ldots,Y_s$, and consider 
${\mathcal{W}}$, the semigroup of words in $Y_1,\ldots,Y_s$. The word $U$ is {\it reduced} if it is  not a power $U_0^k$ for some $k\geq 2$ and shorter word $U_0$. The words $U_1U_2$ and $U_2U_1$ are said to be {\it cyclically equivalent}. Write ${\mathcal{W}}_0$ for the set of cyclic equivalence classes of reduced words. 
Denote by 
$C_{\{(r_1,u_1),\ldots,(r_d,u_d)\}}$ the number of multisets 
$\{U_1,\ldots,U_d;r_1,\ldots,r_d\}$ of cyclic equivalence classes of reduced words (so $U_1,\ldots,U_d$ are distinct elements of ${\mathcal{W}}_0$, with multiplicities $r_1,\ldots,r_d$) such that the specialization $Y_i\mapsto x(i)$ 
($i=1,\ldots,m$) sends $U_j$ to $u_j$ ($j=1,\ldots,d)$; 
note that non-commuting reduced words from different cyclic equivalence classes  may specialize to the same commuting monomial, so here 
$\{(r_1,u_1),\ldots,(r_d,u_d)\}$ is a multiset of pairs 
(i.e. $(r_i,u_i)$ and $(r_j,u_j)$ may be equal for $i\neq j$), 
where $u_i\in\monoms_s$, the set of commuting monomials in $x(1),\ldots,x(s)$.  
By \cite{amitsur} we have 
\begin{equation}\label{eq:amitsur} 
\sigma_l^{\alpha}=\sum(-1)^{l+r_1+\cdots+r_d}
C_{\{(r_1,u_1),\ldots,(r_d,u_d)\}}\sigma_{r_1}(u_1)\cdots\sigma_{r_d}(u_d),
\end{equation} 
where the summation ranges over multisets of pairs 
$\{(r_1,u_1),\ldots,(r_d,u_d)\}$, where $r_i\in\mn$, $u_i\in\monoms_s$ with  
multidegree $\mdeg(u_1^{r_1}\cdots u_d^{r_d})=\alpha$. 
Extending the definition of $\sigma$ by setting $\sigma_l=0$ for $l>n$,  
the equality (\ref{eq:amitsur}) remains valid also if $l>n$. 
For a multidegree $\alpha=(\alpha_1,\ldots,\alpha_s)$ and 
$w_1,\ldots,w_s\in\monoms$ define 
$$S^{\alpha}(w_1,\ldots,w_s)\in \field[e_r(w)\mid w\in\monoms,r=1,\ldots,n]$$ 
as the element obtained by making first the substitution 
$x(i)\mapsto w_i$, $i=1,\ldots,s$ in the right hand side of 
(\ref{eq:amitsur}), and replacing the symbols $\sigma$ by $e$ everywhere. 
Since $\sigma_l$ is identically zero for $l>n$, the equality 
(\ref{eq:amitsur}) shows that 
\begin{equation}\label{eq:sigmarel}
S^{\alpha}(w_1,\ldots,w_s)\in\ker(\phi) 
\quad \mbox{ if }\sum\alpha_i>n.
\end{equation}  

There is another type of relations we shall need. 
Simplify the notation by writing $x=x(1)$ in the special case $m=1$. 
It is well known that  
$\field[V]^{S_n}$ is generated by the algebraically independent elements $\sigma_r(x)$, $r=1,\ldots,n$, hence for each $k\geq 2$ and $j\in\{1,\ldots,n\}$ 
there is a unique $n$-variable polynomial $f_{j,k}$ with integer coefficients 
such that $\sigma_j(x^k)=f_{j,k}(\sigma_1(x),\ldots,\sigma_n(x))$. 
Now given a monomial $w\in\monoms$, consider the element 
$$Q_{j,k}(w)=e_j(w^k)-f_{j,k}(e_1(w),\ldots,e_n(w))\in\ker(\phi).$$ 
For example, it is an easy exercise to verify that 
$$Q_{r,2}(x)=e_r(x^2)-\sum_{i=\max\{0,2r-n\}}^{\min\{2r,n\}}
(-1)^{r+i}e_i(x)e_{2r-i}(x).$$ 
As a first application of the relations (\ref{eq:sigmarel})  
we derive a short proof of the degree bound of Fleischmann \cite{fleischmann} for the generators of $\field[V^m]^{S_n}$. 
We call an element of $\field[V^m]^{S_n}$ {\it indecomposable} if it is not contained in the $\field$-subalgebra generated by strictly lower degree elements, and call it {\it decomposable} otherwise. Obviously, a homogeneous element of $\field[V^m]^{S_n}$ is indecomposable if and only if it is not contained in the square 
$\aug^2$ of the ideal $\aug$ of $\field[V^m]^{S_n}$ spanned by the homogeneous components of positive degree. 

\begin{lemma}\label{lemma:intpart} 
Let $w=x(1)^{\alpha_1}\cdots x(m)^{\alpha_m}\in\monoms$ be a monomial such that 
$\alpha_j\geq n/r$ for some $j$, and $\deg(w)\geq 1+n/r$. 
Then $\sigma_r(w)$ is decomposable.  
\end{lemma} 
\begin{proof} 
Apply induction on $r$. 
Clearly it is sufficient to show that if $kr\geq n$, then 
$\sigma_r(x^ky)\in \aug^2$ (here we write $x,y$ instead of $x(1),x(2)$). 
Indeed, the general case follows on substituting $x,y$ by appropriate monomials. 
From the relation 
$\phi(S^{(kr,r)}(x,y))=0$ we get 
$$(-1)^r\sigma_r(x^ky)+\sum(-1)^{r/i}C_{\{x^{ki}y^i\}}
\sigma_{r/i}(x^{ki}y^i)\in\aug^2,$$ 
where the sum ranges over divisors $i>1$ of $r$ 
(and $C_{\{x^{ki}y^i\}}$ is the number of cyclic equivalence classes of reduced words that specialize to $x^{ki}y^i$). In particular, this sum is empty if $r=1$ and so $\sigma_1(x^ky)\in\aug^2$ for $k\geq n$. 
If $r>1$, then the assumptions in our lemma apply to all summands, hence they  belong to $\aug^2$ by the induction hypothesis, forcing $\sigma_r(x^ky)\in\aug^2$. 
\end{proof} 

\begin{corollary}\label{cor:bound} 
The $\field$-algebra $\field[V^m]^{S_n}$ is generated by the elements 
$\sigma_n(x(i))$, $i\in\{1,\ldots,m\}$, and   
$\sigma_r(x(1)^{\alpha_1}\cdots x(m)^{\alpha_m})$, where  
$r\in\{1,\ldots,n\}$, $\alpha_j<n/r$ for all $j=1,\ldots,m$, and the greatest 
common divisor of $\alpha_1,\ldots,\alpha_m$ is $1$. 
\end{corollary} 

\begin{remark}\label{rem:fleischmannn} {\rm 
In particular, $\field[V^m]^{S_n}$ is generated in degree 
$\leq \max\{m(n-1),n\}$; this is the main result of Fleischmann \cite{fleischmann}, where it is also shown that this bound can not be improved in general. 
Using the method of \cite{fleischmann} we showed in 
\cite{domokos} that if $\field$ is a field of characteristic $p$ and $n=p^k$ 
with $k\in\mn$, then $\sigma_{n-1}(x(1)\cdots x(m))$ is indecomposable. 
Vaccarino \cite{vaccarino} shows that the $\sigma_r(w)$ with $w$ reduced 
(not a power of a monomial with lower degree) generate the ring of multisymmetric functions, and refering to the bound from \cite{fleischmann} 
deduces a finite system of generators. Our Corollary~\ref{cor:bound} reduces further this generating system. A minimal generating system of $\field[V^m]^{S_n}$ appears in \cite{rydh}. }
\end{remark}

\begin{theorem}\label{thm:sigmarelations} 
The ideal $\ker(\phi)$ is generated by 
$$S^{\alpha}(w_1,\ldots,w_s)\qquad \mbox{and}\qquad Q_{r,k}(w), $$
where $s\geq 2$, $w,w_1,\ldots,w_s\in\monoms$;  
$\sum_{i=1}^s\alpha_i>n$, $\sum_{i\neq j}\alpha_i\leq n$ for all $j=1,\ldots,s$;  $r\in\{1,\ldots,n\}$; $k<1+n/r$. 
\end{theorem}  
\begin{proof} Call the {\it weight} of a 
product $\sigma_{r_1}(w_1)\cdots\sigma_{r_s}(w_s)$ the integer 
$r_1+\cdots+r_s$. First we show that if the weight is $>n$, 
then modulo the first type relations in our statement, this product can be rewritten as a linear combination of products with strictly smaller weight. 
It is sufficient to deal with the case when 
$\sum_{i\neq j}r_i\leq n$ for all $j=1,\ldots,s$. 
Set $\alpha=(r_1,\ldots,r_s)$. 
Then $S^{\alpha}(x(1),\ldots,x(s))$ contains the term $e_{r_1}(x(1))\cdots e_{r_s}(x(s))$ with coefficient $\pm 1$. 
In any other term 
$e_{i_1}(u_1)\cdots e_{i_d}(u_d)$ 
of $S^{\alpha}$, at least one of the 
$u_t$ has degree $>1$, hence $i_1+\cdots+i_d<r_1+\cdots+r_s$. 
So the relation $\phi(S^{\alpha}(w_1,\ldots,w_s))=0$ does what we need. 

Next we claim that using our relations, any product 
$\sigma_{r_1}(w_1)\cdots\sigma_{r_d}(w_d)$ with weight $\leq n$ can be 
rewritten as a linear combination of  
$\sigma_{i_1}(u_1)\cdots\sigma_{i_t}(u_t)$ with weight $\leq \sum r_i$ and 
$u_1,\ldots,u_t$ distinct. By Proposition~\ref{prop:basis2} our 
theorem will follow. 
Obviously, it is sufficient to verify the latter claim in the special case 
when all the $w_i$ are powers of the same monomial $w$. 
We may also assume that $w=x$, so we are working in the special case $m=1$. 
An inspection of the proof of Proposition~\ref{prop:basis2} shows 
that $\sigma_{r_1}(x^{j_1})\cdots\sigma_{r_d}(x^{j_d})$ is a linear combination of monomial multisymmetric functions of height $\leq\sum r_i$, and in the $\field$-space spanned by such monomial multisymmetric functions, the products 
$\sigma_{i_1}(x^{k_1})\cdots \sigma_{i_t}(x^{k_t})$ with 
weight $\leq \sum r_i$ and $k_1,\ldots,k_t$ distinct form a basis. 
So we conclude that there are relations among the $\sigma_r(x^k)$ of the desired form. Therefore it suffices to prove the special case $m=1$ of our theorem. 

Since in the special case $m=1$, $\field[V]^{S_n}$ is a polynomial ring generated by $\sigma_1(x),\ldots,\sigma_n(x)$, 
all we need to show is that our relations are sufficient to 
express all $\sigma_r(x^k)$ in terms of the $n$ algebraically independent  generators. 
We apply a double induction: the first induction goes on $\deg(\sigma_r(x^k))=rk$, and the second induction goes on $r$. 
If $(k-1)r<n$, then we may use our relation $\phi(Q_{r,k}(x))=0$ from our statement and we are done.  
Assume $(k-1)r\geq n$. Then we slightly adjust the argument in the proof of 
Lemma~\ref{lemma:intpart}. Write $l$ for the lower integral part of $n/r$. 
Then we have $k>l$ and we have the relation 
$\phi(S^{(lr,r)}(x,x^{k-l}))=0$ in our statement. This relation can be used 
to express 
$$(-1)^r\sigma_r(x^k)+\sum(-1)^{r/i}C_{\{x^{li}y^i\}}
\sigma_{r/i}(x^{ki})$$ 
(where the sum ranges over divisors $i>1$ of $r$) by 
lower degree elements. The statement then follows by the two 
induction hypotheses. 
\end{proof}

\begin{remark}\label{rem:ideal} {\rm 
Although the above theorem is far from yielding a finite presentation of 
$\field[V^m]^{S_n}$, at least it gives a uniform presentation in the sense 
that all relations are obtained from finitely many types by substituting different monomials (note that the restrictions made on $\alpha$ 
force $\sum_{i=1}^s \alpha_i\leq 2n$ and $s\leq n+1$). 
For comparison we mention that in \cite{vaccarino} the author considers 
the generating system $\sigma_r(w)$ where $w$ is reduced 
(i.e. is not a power of a monomial with smaller degree), and describes the relations among these generators in the following sense: he works in the ring of multisymmetric functions in infinitely many variables, where the analogues of 
$\sigma_r(w)$ ($r\in\mn$, $w\in\monoms$ reduced) form an algebraically independent generating system, and observes that the unique expressions of the monomial multisymmetric functions with height $>n$ in terms of these generators form a $\field$-basis in the space of relations among these generators. 
Implementing this procedure in practice one needs to produce formulae like our $Q_{r,k}(w)$ but with no bound on $r$ and $k$, and 
one needs to produce expressions similar to $S^{\alpha}(w_1,\ldots,w_s)$ with 
no bound on $\alpha$. } 
\end{remark} 


\section{Calculations}\label{sec:calculations} 

Throughout this section we assume that $n!$ is invertible in 
$\field$, so $\field[V^m]^{S_n}$ is Cohen-Macaulay 
(see for example \cite{dk} as a general reference for the invariant theory of finite groups). 
To illustrate how to proceed to find a minimal presentation we do calculations in concrete examples. 
An obvious homogeneous system of parameters {\it (primary generators)} in $\field[V^m]^{S_n}$ is 
$$P=\{[x(i)],[x(i)^2],\ldots,[x(i)^n]\mid i=1,\ldots,m\}.$$ 
Write $\parameterid$ for the ideal of $\field[V^m]^{S_n}$ generated by $P$, 
and write $\field[P]$ for the polynomial subalgebra of $\field[V^m]^{S_n}$ 
generated by $P$.  
Note that to find {\it secondary generators} we need to get a vector space 
basis modulo $\parameterid$ in $\field[V^m]^{S_n}$. 

\begin{lemma}\label{lemma:longtrace} 
Assume that $n!$ is invertible in $\field$. 
Let $w$ be a monomial having degree $\geq n$ in one of the variables 
$x(1),\ldots,x(m)$, and having total degree $\geq n+1$. 
Then $[w]$ belongs to $\parameterid$. 
\end{lemma} 

\begin{proof} Assume for example that $w$ has degree $\geq n$ in $x(1)$. 
Then $w$ can be written as a product of $n+1$ factors as 
$x(1)x(1)\cdots x(1)u$ for some non-empty monomial $u$. 
The relation 
$\Psi_{n+1}(x(1),x(1),\ldots,x(1),u)=0$ 
shows the claim, since each nontrivial partition of the multiset 
$\{x(1),\ldots,x(1),u\}$ contains a part consisting solely of $x(1)$s. 
\end{proof}

\begin{remark}\label{rem:seondarybound} 
{\rm From Lemma~\ref{lemma:longtrace} and Proposition~\ref{prop:basis} 
we immediately get the upper bound $n(n-1)m$ for the degrees of the 
secondary generators. However, the general method of Broer 
(see Theorem 3.9.8 in \cite{dk}) 
yields the better bound $\frac{1}{2}n(n-1)m$.}
\end{remark}  

\subsection{The case $n=2$}\label{subsec:n=2} 

In this subsection we assume $\char(\field)\neq 2$. 
The fundamental multilinear trace identity for diagonal $2\times 2$ matrices is 
\begin{equation}\label{eq:fundtraceid2} 
[xyz]=\frac{1}{2}([xy][z]+[xz][y]+[yz][x]-[x][y][z])
\end{equation}
Substitute $z\mapsto zw$ in (\ref{eq:fundtraceid2}), 
and eliminate traces of degree $3$ using (\ref{eq:fundtraceid2}) 
to get 
\begin{eqnarray}\label{eq:2trace4} 
[xyzw]=
\frac{1}{4}([xz][y][w]+[xw][y][z]
+[yz][x][w]+[yw][x][z])
\notag \\
+\frac{1}{2}[xy][zw]-\frac{1}{2}[x][y][z][w]
\end{eqnarray} 
Exchanging the variables $y$ and $z$ in (\ref{eq:2trace4}) we obtain 
\begin{eqnarray}\label{eq:2trace4-1}
[xzyw]=
\frac{1}{4}([xy][z][w]+[xw][y][z]
+[yz][x][w]+[zw][x][y])
\notag \\
+\frac{1}{2}[xz][yw]-\frac{1}{2}[x][y][z][w]
\end{eqnarray} 
The difference of  (\ref{eq:2trace4}) and (\ref{eq:2trace4-1}) yields 
\begin{equation}\label{eq:2basic} 
[xz][yw]-[xy][zw]=\frac{1}{2}([xy][z][w]+[zw][x][y]-[xz][y][w]-[yw][x][z])
\end{equation} 
(compare with (4.17) in \cite{olver}). 
The specialization $z=w$ in (\ref{eq:2basic}) leads to 
\begin{equation}\label{eq:2square}
[xz][yz]=[xy][zz]+\frac{1}{2}([xy][z][z]+[zz][x][y]-[xz][y][z]-[yz][x][z])
\end{equation}
The $\mn_0^m$-graded Hilbert series of $\field[V^m]^{S_2}$ can be written as a rational function by 
Molien's formula: 
$$H(\field[V^m]^{S_2};t_1,\ldots,t_m)
=\frac{\sum_{i=0}^{[\frac{m}{2}]}e_{2i}(t_1,\ldots,t_m)}
{\prod_{j=1}^m(1-t_j)(1-t_j^2)},$$
where $e_i(t_1,\ldots,t_m)$ is the $i$th elementary symmetric function in 
$t_1,\ldots,t_m$. 
It is easy to describe the Hironaka decomposition of $\field[V^m]^{S_2}$. 
The secondary generators are 
$$S=\{[x(i_1)x(i_2)]\cdots[x(i_{2k-1})x(i_{2k})]\mid 
1\leq i_1<\cdots<i_{2k}\leq m\}.$$
The substitution $x\mapsto x(i)$, $y\mapsto x(j)$, $z\mapsto x(k)$ 
in (\ref{eq:2square}) 
gives a relation showing that $[x(i)x(k)][x(j)x(k)]$ 
is contained in $\parameterid$ for all $i,j,k$. 
Substitutions $\{x,y,z,w\}\to \{x({j_1}),\ldots,x(j_{2k})\}$ 
in the relation (\ref{eq:2basic}) show that 
$$[x(j_1)x(j_2)]\cdots[x(j_{2k-1})x(j_{2k})]-
[x(j_{\pi(1)})x(j_{\pi(2)})]\cdots [x(j_{\pi(2k-1)})x(j_{\pi(2k)})]$$ 
belongs to 
$\langle P\rangle$  
for an arbitrary permutation $\pi$. 
These relations imply that the algebra of multisymmetric functions is generated by $S$ as a module over the polynomial ring $\field[P]$. 
The Hilbert series shows that it is a free module. 
These considerations show also that the specializations 
$\{x,y,z,w\}\to \{x(1),\ldots,x(m)\}$ in (\ref{eq:2basic}) generate the ideal of 
relations among the generators $[x(i)],[x(k)x(l)]$, and that a minimal system 
of relations consists of relations of degree $4$. 

\subsection{The case $n=3$, $m=2$}\label{subsec:n=3} 

In this subsection we assume $\char(\field)\neq 2,3$. 
The fundamental identity $\Psi_4(x,y,z,w)=0$ 
takes the form 
\begin{eqnarray*}
[xyzw]=\frac{1}{3}([xyz][w]+[xyw][z]+[xzw][y]+[yzw][x])\\
+\frac{1}{6}([xy][zw]+[xz][yw]+[xw][yz]+[x][y][z][w])
\\-\frac{1}{6}
([xy][z][w]+[xz][y][w]+[xw][y][z]+
\\ +[yz][x][w]+[yw][x][z]+[zw][x][y])
\end{eqnarray*}
Consider the following consequences: 
\begin{eqnarray} \label{eq:3rel}
& \Psi_4(x,x,x,x)=0 \qquad \Psi_4(y,y,y,y)=0
\qquad \Psi_4(x,x,x,y)=0 
\\ \notag & \Psi_4(y,y,y,x)=0 \qquad \Psi_4(x,x,x,xy)=0 \qquad \Psi_4(y,y,y,xy)=0  
\\ \notag & \Psi_4(x,x,y,y)=0 \qquad \Psi_4(x,x,x,y^2)=0 \qquad \Psi_4(y,y,y,x^2)=0
\\ \notag & \Psi_4(x,x,x,y^3)=0 \qquad \Psi_4(x,x,x,xy^2)=0 \qquad \Psi_4(y,y,y,yx^2)=0 
\end{eqnarray} 
and 
\begin{eqnarray} \label{eq:3rel2} 
& \Psi_4(x^2,x,y,y)=0 \qquad \Psi_4(y^2,y,x,x)=0
\\ \notag &\Psi_4(x^2,x^2,y,y)=0 \qquad \Psi_4(y^2,y^2,x,x)=0 
\\ \notag & \Psi_4(x^2,x,y^2,y)=0
\end{eqnarray} 
Let us use the notation $f\equiv g$ to indicate that $f,g\in\field[V^2]^{S_3}$  
are congruent modulo $\parameterid$. 
The above 17 relations in the order of listing imply the following: 
\begin{eqnarray*} \label{eq:3congruences} 
&[x^4]\equiv 0 \qquad [y^4]\equiv 0 
\qquad [x^3y]\equiv 0 
\\  & [xy^3]\equiv 0  \qquad [x^4y]\equiv 0 \qquad [xy^4]\equiv 0 
\\  & [x^2y^2] \equiv \frac{1}{3} [xy][xy] 
\qquad [x^3y^2]\equiv 0 \qquad [y^3x^2]\equiv 0
\\  & [x^3y^3]\equiv 0 \qquad [x^4y^2]\equiv 0\qquad [y^4x^2]\equiv 0
\end{eqnarray*} 
and 
\begin{eqnarray*}
& [x^3y^2]\equiv \frac{1}{3}[x^2y][xy] \qquad 
 [y^3x^2]\equiv \frac{1}{3}[y^2x][yx] 
\\ & [x^4y^2]\equiv \frac{1}{3}[x^2y][x^2y] 
\qquad [x^2y^4]\equiv \frac{1}{3}[xy^2][xy^2] 
\\ & [x^3y^3]\equiv \frac{1}{6}[x^2y][xy^2]+
\frac{1}{6}[x^2y^2][xy]
\end{eqnarray*}
Recall that $\field[V^2]^{S_3}$ is minimally generated as an algebra by 
$P\cup\{[xy],[x^2y],[xy^2]\}$. 
The above congruences obviously show that as a module over 
$\field[P]$, the algebra $\field[V^2]^{S_3}$ is spanned by  
$$S=\{1,[xy],[x^2y],[xy^2],[xy]^2,[x^2y][xy^2]\}.$$ 
Using Molien's formula one can easily compute 
$$H(\field[V^2]^{S_3};t_1,t_2)=
\frac{1+t_1t_2+t_1^2t_2+t_1t_2^2+t_1^2t_2^2+t_1^3t_2^3}
{(1-t_1)(1-t_1^2)(1-t_1^3)(1-t_2)(1-t_2^2)(1-t_2^3)},$$ 
and this implies that $\field[V^2]^{S_3}$ is a free $\field[P]$-module 
generated by $S$.  
The above considerations imply that the ideal of 
relations among the 21 generators 
$$\{[x^iy^j]\mid 0\leq i,j\leq 4, \quad i+j\leq 6,\quad (i,j)\neq (0,0)\}$$ 
is generated by the 17 relations (\ref{eq:3rel}) and (\ref{eq:3rel2}). 
The relations (\ref{eq:3rel}) can be used to eliminate 
the $12$ superfluous generators of degree $\geq 4$ in the relations (\ref{eq:3rel2}). 
This way we get 5 defining relations for the minimal generating system 
$\{[x^iy^j]\mid i+j\leq 3\}$. 
The bidegrees of these relations are 
$(3,2)$, $(2,3)$, $(4,2)$, $(2,4)$, $(3,3)$. 
It is easy to see that these 5 relations minimally generate the ideal of relations. 
(In subsequent work, the minimal presentation of $\mc[V^m]^{S_3}$ is determined in \cite{domokos-puskas} for arbitrary $m$.)


\end{document}